\documentclass{aip-cp}

\usepackage{bm}
\usepackage[numbers]{natbib}
\usepackage{rotating}
\usepackage[utf8]{inputenc}

\begin{document}

\newcommand{\estim}[2]{\ensuremath{\left\{#1\right\}_{#2}}}

\title{Plausible Inference From The Idea Of Estimation}

\author[aff1]{Sergio Davis\corref{cor1}}

\affil[aff1]{Comisión Chilena de Energía Nuclear, Casilla 188-D, Santiago, Chile.}
\corresp[cor1]{Corresponding author: sdavis@cchen.cl}

\maketitle

\begin{abstract}
The probability axioms by R. T. Cox can be regarded as the modern foundations of Bayesian inference, 
the idea of assigning degrees of belief to logical propositions in a manner consistent with Boolean logic.
In this work it is shown that one can start from an alternative point of view, postulating the 
existence of an operation for \textit{estimating magnitudes given a state of knowledge}. The properties 
of this operation can be established through basic consistency requirements. Application of this 
estimation operation to the truth values of propositions in an integer representation directly yields 
non-negative, normalized quantities following the sum and product rules of probability.
\end{abstract}

\section{INTRODUCTION}

The idea of plausible inference was masterfully cast into a formal theory by Richard T. Cox in 1946~\cite{Cox1946}. In 
his work, Cox postulates requirements (\emph{desiderata}) for an extension of propositional logic, in the form of a 
mapping between logical propositions and non-negative real numbers, thus introducing a unique gradation between absolute 
falsehood and absolute truth. This algebra of plausible inference coincides with the rules for the manipulation of 
conditional probabilities in the frequentist approach, and therefore lends support to the Bayesian interpretation of 
probabilities as degrees of belief.

Despite the clarity of Cox's work, and later developments by Jaynes~\cite{Jaynes2003} among others~\cite{Sivia2006} to
promote the conceptual power and practical utility of Bayesian ideas, their adoption has not been straightforward, and 
the subject of probabilities still seems clouded in mystery. See for instance the work by Cousins~\cite{Cousins1994} 
on the adoption of Bayesian methods by physicists. In contrast with the standard academic formation in probability and statistics, 
our intuition in Science and particularly in Physics deals with quantities -- properties of a system such as mass or 
temperature, or parameters in a phenomenological model-- to be estimated from a set of measurements. In this case probability 
only enters indirectly, as the ``weight'' we postulate for every admissible value of the quantity. Based on these 
ideas we could imagine a theory of estimation of numerical quantities, which does not directly assign degrees of belief 
to propositions but rather assigns an \emph{estimated value} to a quantity given some information. The advantage of a 
formulation centered on estimation is that it has the potential to bypass the philosophical issues attached to the 
concept of probability.

In this work, the starting elements for such a theory of estimation are discussed, in the form of requirements of 
consistency. In the proposed framework it is possible to speak of the estimation of the value of quantity -- either a 
random variable or an unknown parameter -- in agreement with Bayesian ideas but without referring to an underlying
probability. Probabilities as known in the standard Bayesian formulation emerge when this theory is applied to the 
estimation of truth values in their integer representation (binary digits). For both discrete and continuous variables 
the estimation operation is shown to be equivalent to the usual expectation, under a probability defined as the estimation 
of a delta function in the corresponding space.

\section{CONSISTENCY REQUIREMENTS FOR ESTIMATION}

We will introduce the notation $\estim{x}{I}$ to denote the \emph{estimated value} of the quantity $x$ under the state
of knowledge $I$. This notation immediately tells us that the notion of an absolute estimation of $x$ is meaningless; 
every estimation needs some input. The quantity $x$ to be estimated may have discrete or continuous values, but its 
estimation $\estim{x}{I}$ will be in general a real number, in order to allow estimations to be closer to one particular 
discrete value rather than others in a gradual manner. 

The value $\estim{x}{I}$ only depends on the parameters (if any) defined by the proposition $I$, and any symbol inside 
the curly braces not defined elsewhere (such as $x$ in $\estim{x}{I}$) is assumed to be an unknown quantity, subject to 
estimation. 
 
For instance, if $x$, $y$ and $w$ are unknowns, 

\begin{enumerate}
\item[(a)] the expression $\estim{w}{I}$ is just a number, 
\item[(b)] the expression $\estim{x\cdot y}{x=x_0}$ is a function of $x_0$ only, and 
\item[(c)] the expression $\estim{f(x,y)}{y>y_0}$ is a function of $y_0$.
\end{enumerate}
In general, we can understand the estimation as the incomplete evaluation of a function, in the 
case when its arguments are not fully known.

Now we will postulate some consistency requirements for the estimation operation, in order for it to correspond with
common sense, in the same spirit as Cox did with probabilities. In the statement of all the requirements below, $x$, $y$, 
$z$ and $w$ are unknown quantities, $\alpha$ is a known constant with value $\alpha_0$ and the state of knowledge $I$ 
is considered to be arbitrary.

\begin{enumerate}

\item[(0)] Estimation of a fully known quantity yields that quantity. That is, if $f$ is a function of $x$,

\begin{equation}
\estim{f(x)}{x=x_0, I} = f(x_0).
\end{equation}
This is the fundamental property required of an estimation: it has to agree with simple evaluation of a function 
when all its arguments are known. A particular implication of this requirement is that the estimated value of 
a known constant $\alpha$ is the constant itself, 

\begin{equation}
\estim{\alpha}{\alpha=\alpha_0,I}=\alpha_0.
\end{equation}

\vspace{12pt}

\item[(1)] The estimated value of a quantity cannot lie outside its admissible interval. As noted before, the 
estimation of a discrete value does not necessarily correspond to one of its allowed values. However, it must 
lie inside the interval between its minimum and maximum discrete values. That is,

\begin{equation}
x \in \{x_1, \ldots, x_N\} \Longrightarrow \min(x_1, \ldots, x_N) \leq \estim{x}{I} \leq \max(x_1, \ldots, x_N).
\end{equation}

For a continuous variable $x$ with values inside an interval $\Omega$, the requirement is that $x \in \Omega 
\Longrightarrow \estim{x}{I} \in \Omega$.

\vspace{12pt}

\item[(2)] Estimation is linear, so that the estimation of a quantity scaled by a known constant is just the 
scaled estimation of the ``bare'' quantity,

\begin{equation}
\estim{\alpha x}{\alpha=\alpha_0,I} = \alpha_0\estim{x}{\alpha=\alpha_0,I}, \nonumber \\
\end{equation}
and the estimation of a sum is the sum of estimations (under the same state of knowledge),

\begin{equation}
\estim{x + y}{I} = \estim{x}{I}+\estim{y}{I}.
\end{equation}

We require the estimation to be linear because one must be able to estimate physical quantities 
regardless of the chosen measurement unit. This implies that a scaling of the quantity can only rescale
the estimation, and a shift of the origin can only shift the resulting estimation. It also must be the case 
that estimating two quantities separately and adding them is equivalent to estimating their sum simultaneously.

\vspace{12pt}

\item[(3)] Partial estimation removes free parameters. \\

Let us define a partial estimation $\estim{y}{x,I}$ as the result of estimating the quantity $y=y(x,z,w)$ 
under the assumption of known value of a parameter $x$ and state of knowledge $I$. This partial estimation is of 
course a function of $x$ only, into which the information about $z$ and $w$ has been incorporated by the 
estimation procedure. We require that the remaining free parameter $x$ can be estimated as well, by 
applying a second estimation operation, and that the result is equivalent to the estimation of $y$ under 
$I$ without any intermediate steps. That is,

\begin{equation}
\estim{y}{I} = \estim{ \estim{y}{x,I} }{I}.
\end{equation}
If estimation is analogous to partial evaluation of a function, then a second estimation operation effectively 
removes all remaining free arguments. In fact, in the right-hand side the quantity $x$ which is a fixed parameter in the 
inner estimation becomes an unknown for the purposes of the outer estimation.

\end{enumerate}

We will see in the next section that these requirements by themselves are sufficient to uniquely define the correct 
implementation of the estimation operation. Moreover, we will show that the sum and product rules of 
probability theory follow when estimation is used on binary values originating from propositional logic.

\section{PROBABILITY THEORY}

It is well-known that the Boolean operations $\wedge$ (logical \texttt{and}), $\vee$ (logical \texttt{or})
and $\neg$ (logical \texttt{not}) can be mapped into arithmetic operations on $\{0, 1\}$. To illustrate this, a convenient 
notation is the mapping $n(A) \in \{0, 1\}$, with $A$ a logical proposition, so that $n($True$)=1$ and $n($False$)=0$. 
Note that

\begin{equation}
n: \{\texttt{False}, \texttt{True}\} \rightarrow \{0, 1\}
\end{equation}
goes strictly from logical propositions to the integers zero and one, not to values in between. Given these definitions, 
it is straightforward to encode the logical negation $\neg$ as

\begin{equation}
n(\neg A) = 1 - n(A),
\label{eq_bit_negation}
\end{equation}
and the logical operators $\vee,\wedge$ (shown in Table \ref{tbl_truth}) as

\begin{eqnarray}
n(A \vee B) = n(A) + n(B) - n(A)n(B), \\
n(A \wedge B) = n(A)n(B).
\label{eq_bit_andor}
\end{eqnarray}
The values of $n(A\vee B)$ and $n(A \wedge B)$ are reproduced for convenience in Table \ref{tbl_bits}. Here the 
substracted term $n(A)n(B)$ in the expression for $n(A \vee B)$ is used to implement sum \emph{modulo} 2.

\begin{table}[h]
\caption{Truth values of the binary Boolean operators $\wedge$ (logical \texttt{and}) and $\vee$ (logical \texttt{or}).}
\label{tbl_truth}
\begin{tabular}{cccc}
\hline
$A$ & $B$ & $A \vee B$ & $A \wedge B$ \\
\hline
False & False & False & False \\
False & True  & True  & False \\
True  & False & True  & False \\
True  & True  & True  & True \\
\hline
\end{tabular}
\end{table}

\begin{table}[h]
\caption{Binary values of the sum \emph{modulo} 2 and product operators.}
\label{tbl_bits}
\begin{tabular}{cccc}
\hline
$n(A)$ & $n(B)$ & $n(A)+n(B)-n(A)n(B)$ & $n(A)n(B)$ \\
\hline
0  & 0 & 0 & 0 \\
0  & 1 & 1 & 0 \\
1  & 0 & 1 & 0 \\
1  & 1 & 1 & 1 \\
\hline
\end{tabular}
\end{table}

Now we proceed to apply the estimation operator $\estim{\cdot}{I}$ to the integer $n(A)$, where $I$ is a proposition 
connected in an arbitrary way to propositions $A$ and $B$. The resulting quantity $\estim{n(A)}{I}$ is guaranteed to be 
a real number between 0 and 1 because of the discrete form of Requirement 1. Therefore, Equations \ref{eq_bit_negation} 
to \ref{eq_bit_andor} imply that the estimations $\estim{n(A)}{I}$, $\estim{n(B)}{I}$ and $\estim{n(A)n(B)}{I}$ are 
such that 

\begin{eqnarray}
\estim{n(\neg A)}{I} = 1 - \estim{n(A)}{I}, \\
\estim{n(A \vee B)}{I} = \estim{n(A)}{I} + \estim{n(B)}{I} - \estim{n(A)n(B)}{I}, \\
\estim{n(A \wedge B)}{I} = \estim{n(A)n(B)}{I}.
\end{eqnarray}
 
At this point it is tempting to read $\estim{n(A)}{I}$ as the Bayesian probability $P(A|I)$. However, the 
product rule is not manifestly reproduced, as is not apparent how $\estim{n(A \wedge B)}{I}=\estim{n(A)n(B)}{I}$ 
would be equivalent to $\estim{n(A)}{I}\estim{n(B)}{A, I}$. To prove that $\estim{n(A)n(B)}{I}$ in fact reduces to 
$\estim{n(A)}{I}\estim{n(B)}{A, I}$, we introduce the auxiliary variables $a=n(A)$ and $b=n(B)$ with $a,b \in \{0, 1\}$. 
Then according to Requirement 3 we see that

\begin{equation}
\estim{a\cdot b}{I}=\estim { \estim{a \cdot b}{a,I} }{I}.
\end{equation}
Now in the right-hand side, $a$ is a constant inside the inner estimation and can be extracted out, so we have

\begin{equation}
\estim { \estim{a \cdot b}{a,I} }{I} = \estim{ a\cdot\estim{b}{a,I}}{I} = \estim{G(a)}{I},
\end{equation}
where $G(a)=a\cdot\estim{b}{a,I}$. Using the fact that $a \in \{0, 1\}$ we can rewrite $G$ as the equivalent expression 

\begin{equation}
G^*(a)=a\cdot\estim{b}{a=1,I},
\end{equation}
so that $G(0)=G^*(0)=0$ and $G(1)=G^*(1)=\estim{b}{a=1,I}$. We have then that

\begin{equation}
\estim{a\cdot b}{I}=\estim{G(a)}{I}=\estim{G^*(a)}{I}=\estim{a\cdot\estim{b}{a=1,I}}{I},
\end{equation}
but because $\estim{b}{a=1,I}$ is a constant, we can move it outside the outer estimation and recover

\begin{equation}
\estim{a\cdot b}{I} = \estim{a}{I}\estim{b}{a=1,I},
\end{equation}
which in the previous notation is

\begin{equation}
\estim{n(A)\cdot n(B)}{I} = \estim{n(A)}{I}\estim{n(B)}{A,I}.
\end{equation}
This proves that the identification $P(A|I)=\estim{n(A)}{I}$ leads to the correct form of the product and sum 
rules of probability,

\begin{eqnarray}
P(\neg A|I) = 1 - P(A|I), \\
P(A \vee B|I) = P(A|I) + P(B|I) - P(A\wedge B|I), \\
P(A \wedge B|I) = P(A|I)P(B|A,I).
\end{eqnarray}
without assigning any meaning to the probabilities themselves. Here probabilities are simply plausible estimations 
of binary digits that follow the arithmetic operations \emph{modulo} 2 associated with Boolean logic, and therefore
reveal themselves as an extension of this logic.

\section{CONNECTION BETWEEN ESTIMATION AND EXPECTATION}

Up to this point we do not know if the estimation operation is unique or what is the mathematical machinery needed to 
perform estimation in practice. Fortunately, both questions are straightforward to answer. In order to do this, consider 
first a discrete variable $x \in \{x_1, x_2, \ldots, x_N\}$ and its estimation $\estim{x}{I}$ given $I$. We can always 
rewrite the value of $x$ by including the complete set of values $\{x_1, x_2, \ldots, x_N\}$,

\begin{equation}
x = \sum_{i=1}^N x_i\delta(x_i, x),
\end{equation}
so that after taking estimation under $I$ on both sides, we have

\begin{equation}
\estim{x}{I} = \estim{\sum_{i=1}^N x_i\delta(x_i, x)}{I}.
\end{equation}
By Requirement 2, the estimation operation is linear so that the sum and $x_i$ can fall outside it,

\begin{equation}
\estim{x}{I} = \estim{\sum_{i=1}^N x_i\delta(x_i, x)}{I} = \sum_{i=1}^N x_i\estim{\delta(x_i, x)}{I} = \sum_{i=1}^N x_i p_i
\end{equation}
and so we see that the estimation is a linear combination of all the admissible values of $x$ weighted by non-negative numbers 
$p_i=\estim{\delta(x_i, x)}{I}$. The weights $p_i$ are non-negative due to Requirement 1: 

\begin{equation}
\delta(x_i, x) \in \{0, 1\} \Longrightarrow 0 \leq \estim{\delta(x_i,x)}{I} \leq 1.
\end{equation}
Moreover, they are properly normalized because 

\begin{equation}
\sum_{i=1}^N p_i = \sum_{i=1}^N\estim{\delta(x_i,x)}{I}=\estim{\sum_{i=1}^N\delta(x_i,x)}{I}=1,
\end{equation}
so in summary, $\estim{x}{I}\sum_i x_i p_i$ is precisely the standard expectation $\left<x\right>_I = \sum_{i=1}^N x_i P(x=x_i|I)$ 
under a probability 

\begin{equation}
P(x=x_i|I)=\estim{\delta(x_i,x)}{I}.
\end{equation}
This probability is consistent with our previous definition $P(A|I)=\estim{n(A)}{I}$ for a logical proposition $A$, because
$\delta(x_i,x)$ is 1 if $x=x_i$ and 0 otherwise. In other words, $\delta(x_i,x)$ is equal to $n(x=x_i)$.
In the case of a continuous variable $x$ we can do a similar expansion,

\begin{equation}
\estim{x}{I} = \estim{\int_a^b dx'\; x'\delta(x'-x)}{I} = \int_a^b dx'\; x'\estim{\delta(x'-x)}{I} = \int_a^b dx'\; x'\;p(x'),
\end{equation}
with $p(x')=\estim{\delta(x'-x)}{I}$ a non-negative number corresponding to the probability density $P(x=x'|I)$. Again, given 
the values of Dirac's delta it follows that $0 \leq \estim{\delta(x'-x)}{I} < \infty$. A general rule seems to emerge, 
associating probabilities to estimations of two-valued functions such as Kronecker's delta, Dirac's delta or the Heaviside step 
function, which encode logical propositions.

\section{CONCLUSIONS}

We have presented an alternative derivation of probability theory, starting from the idea of estimation given information.
In a similar way to Cox's original arguments for his formulation of probability theory, based on consistency with Boolean logic, 
we have shown that a consistent definition of estimation applied to truth values in their integer representation leads to non-negative, 
normalized quantities that operate using the sum and product rule of probability theory. 

There are several ways to extend these results. For instance, we have defined estimation only for discrete and continuous (real) 
numbers, but it is a relatively simple task to generalize these ideas to vectors, matrices, functions among other mathematical objects 
which can be decomposed as linear combinations in a particular basis. This provides a solution to the problem of defining 
an algebra of probabilities on abstract spaces such as dynamical trajectories without describing them in terms of logical propositions.

\section{ACKNOWLEDGMENTS}

SD acknowledges support from FONDECYT grant 1140514.

\bibliographystyle{aipnum-cp}
\bibliography{estimation}

\end{document}